\documentclass[11pt]{article}
\usepackage{tipa}
\usepackage[all,import]{xy}
\usepackage{graphicx,color}
\usepackage{amsmath}
\usepackage{amsbsy}
\usepackage{amssymb}
\usepackage{cases}
\usepackage{enumerate}
\usepackage{bm}
\usepackage[colorlinks,linkcolor=black,anchorcolor= black,citecolor= black,urlcolor=black]{hyperref}
\RequirePackage{hyperref}
\RequirePackage{hypernat}

\usepackage{amsmath}

\usepackage{amsthm}
\usepackage{multirow}
\usepackage{subfigure}
\usepackage{colortbl}
\usepackage{algorithm}
\usepackage{algorithmic}

\usepackage{times}
\usepackage{mathptmx}

\definecolor{mygray}{gray}{.9}

\newcommand{\Gammad}{\text{Gamma}}

\DeclareMathOperator\arctanh{arctanh}
\DeclareMathOperator\sech{sech}

\makeatletter
\newcommand*{\rom}[1]{\expandafter\@slowromancap\romannumeral #1@}
\makeatother

\begin{document}

\theoremstyle{definition}
\newtheorem{assumption}{Assumption}
\newtheorem{theorem}{Theorem}
\newtheorem{lemma}{Lemma}
\newtheorem{example}{Example}
\newtheorem{definition}{Definition}
\newtheorem{corollary}{Corollary}

\def\letas{\mathrel{\mathop{=}\limits^{\triangle}}}
\def\ind{\begin{picture}(9,8)
         \put(0,0){\line(1,0){9}}
         \put(3,0){\line(0,1){8}}
         \put(6,0){\line(0,1){8}}
         \end{picture}
        }
\def\nind{\begin{picture}(9,8)
         \put(0,0){\line(1,0){9}}
         \put(3,0){\line(0,1){8}}
         \put(6,0){\line(0,1){8}}
         \put(1,0){{\it /}}
         \end{picture}
    }

\def\AVar{\text{AsyVar}}
\def\Var{\text{Var}}
\def\Cov{\text{Cov}}
\def\sumn{\sum_{i=1}^n}
\def\summ{\sum_{j=1}^m}
\def\convergeas{\stackrel{a.s.}{\longrightarrow}}
\def\converged{\stackrel{d}{\longrightarrow}}
\def\iidsim{\stackrel{i.i.d.}{\sim}}
\def\indsim{\stackrel{ind}{\sim}}
\def\asim{\stackrel{a}{\sim}}
\def\d{\text{d}}

\setlength{\baselineskip}{2\baselineskip}

\setcounter{MaxMatrixCols}{20}

\title{\bf Three Occurrences of the Hyperbolic-Secant Distribution}
\author{Peng Ding \\ 
Department of Statistics, Harvard University, One Oxford Street, Cambridge 02138 MA\\
Email: pengding@fas.harvard.edu
}
\date{}
\maketitle

\begin{abstract}
Although it is the generator distribution of the sixth natural exponential family with quadratic variance function, the Hyperbolic-Secant distribution is much less known than other distributions in the exponential families.
Its lack of familiarity is due to its isolation from many widely-used statistical models. We fill in the gap by showing
three examples naturally generating the Hyperbolic-Secant distribution, including
Fisher's analysis of similarity between twins, the Jeffreys' prior for contingency tables, and invalid instrumental variables.

\end{abstract}

\section{The Hyperbolic-Secant Distribution: A Review}

\begin{figure}
\centering
\includegraphics[width = 0.5\textwidth]{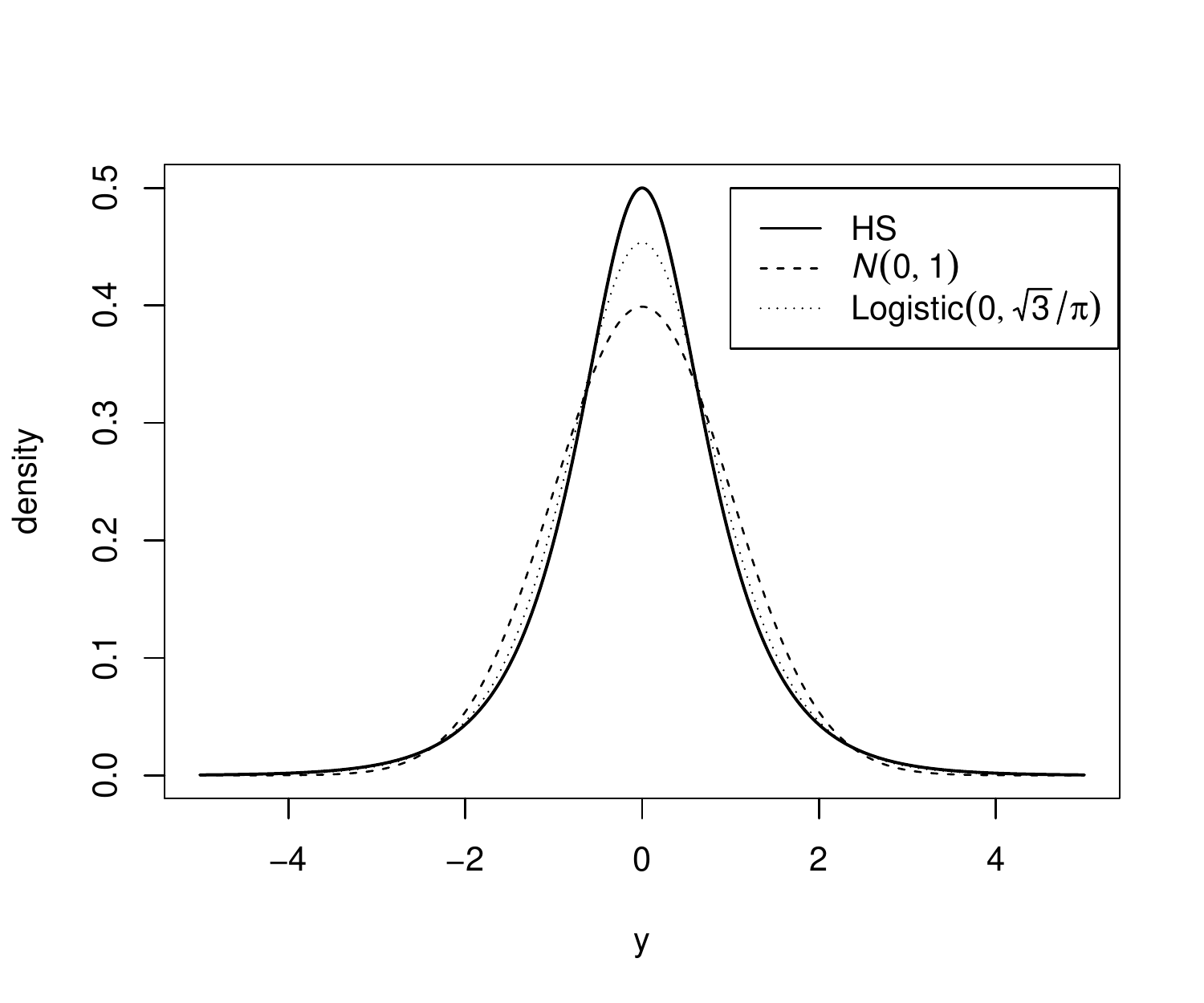}
\caption{HS, Normal and Logistic with mean $0$ and variance $1$}\label{fg::hs_normal}
\end{figure}

The Hyperbolic-Secant (HS) distribution is a bell-shaped distribution with mean $0$ and variance $1$. It can be represented through the standard Cauchy distribution as 
\begin{eqnarray}\label{eq::HS-Cauchy}
Y\sim \frac{2}{\pi} \log |C|,
\end{eqnarray}
where $C$ has the standard Cauchy distribution, and $Y$ has the HS distribution. The density of $Y$ is
$$
f_Y(y) = \frac{1}{2} \text{sech}\left(  \frac{\pi y}{2} \right) = \frac{1}{e^{\pi y/2} + e^{-\pi y/2}}, \quad  y\in (-\infty, +\infty),
$$
the moment generating function of $Y$ is $\sec(t)$, and the characteristic function of $Y$ is $\sech(t)$, for $|t|<\pi/2.$
Figure \ref{fg::hs_normal} compares three densities: HS, $N(0,1)$, and Logistic$(0,\sqrt{3}/\pi)$, all of which have mean $0$ and variance $1$.
The density of the HS distribution has a sharper peak near $0$ and heavier tails than the density of $N(0,1)$, and the density of Logistic$(0,\sqrt{3}/\pi)$ lies between them. Throughout the paper, $Y$, $C$ and $Z$ will be used exclusively to denote random variables with the HS distribution, the standard Cauchy distribution, and the standard Normal distribution, respectively.

The HS distribution is the generator distribution of the sixth natural exponential family (NEF) with quadratic variance function (QVF) (Morris 1982; Morris and Lock 2009).  
Morris (1982) first gave the name NEF-QVF for natural exponential families with the variances as quadratic functions of the mean parameters. The other five NEF-QVFs are Normal, Poisson, Gamma, Binomial, and Negative Binomial, which are much more widely-used and well-known than the HS distribution.

Perks (1932) first derived a family of generalized HS distributions, which could better fit the observed data of the rate of mortality in actuarial science.
Talacko (1956) made a connection between Brownian Motion and the HS distribution.
Harkness and Harkness (1968) calculated the moments and cumulants, and discussed statistical inference for a class of generalized HS distributions.
Manoukian and Nadeau (1988) obtained the cumulative distribution function of the sample mean of the HS distribution via its characteristic function.
Vaughan (2002) showed two empirical studies in which models using the HS distribution could fit the tails better than models using Normal distributions.

The fact that the HS distribution is still mysterious to many people, even within the statistics community, is 
partly due to its
lack of connections to other commonly-used statistical models. 
We will show three examples where the HS distribution arises naturally, including Fisher's analysis of similarity between twins, the Jeffreys' prior for contingency tables, and invalid instrument variables.

\section{How Similar are Twins?}

Assume $(X_1,X_2)$ are characteristics of a pair of twins, distributed as
$$
\begin{pmatrix}
X_1  \\
X_2
\end{pmatrix}
\sim
\bm{ \mathcal{N} }_2
\left\{ \begin{pmatrix}
\mu \\ \mu 
\end{pmatrix} ,
\begin{pmatrix}
\sigma^2 & \rho \sigma^2 \\
\rho \sigma^2 & \sigma^2
\end{pmatrix}
 \right\},
$$
where $X_1$ and $X_2$ are symmetric with Pearson correlation coefficient $\rho.$
Thorndike (1905) used the ``intraclass correlation coefficient'' as a measure of similarity between twins, defined as
$$
R = \frac{2(X_1-\mu)(X_2-\mu)}{  (X_1-\mu)^2 + (X_2-\mu)^2 }.
$$
We can show that $-1\leq R \le 1$, and $\pm 1$ are attainable when $X_1$ and $X_2$ are perfectly correlated with $\rho = \pm 1.$
The ``intraclass correlation coefficient'' is an association measure for paired observations, that estimates $\rho$ (when $\mu$ is given). Thorndike (1905) used it to estimate the correlation between individual pairs of twins. 

In the following, we derive the exact distribution of $R$, and show its relationship with $\rho$. We first apply Fisher's z-transformation to $R$, and then show that $\arctanh(R)$ is a location and scale transformation of the HS distribution.
To be more specific, we have
\begin{eqnarray}\label{eq::fisher-z}
\arctanh(R) = \frac{1}{2} \log \frac{1+R}{1-R} = \frac{1}{2} \log \frac{   \{   (X_1 -\mu) + (X_2-\mu)    \}^2    }{   \{  (X_1-\mu) - (X_2-\mu)   \}^2  }
= \log\Big|   \frac{      (X_1 -\mu) + (X_2-\mu)      }{     (X_1-\mu) - (X_2-\mu)  }     \Big|
. 
\end{eqnarray} 
In the above Equation (\ref{eq::fisher-z}), we have $ (X_1-\mu) + (X_2-\mu)  \sim \mathcal{N}\left\{  0, 2\sigma^2 (1+\rho) \right\} \sim \sqrt{ 2\sigma^2 (1+\rho)   }Z_1$, and $ 
 (X_1-\mu) - (X_2-\mu) \sim \mathcal{N}\left\{ 0, 2\sigma^2(1-\rho)\right\}  \sim \sqrt{2\sigma^2 (1 - \rho)} Z_2$, where $Z_i$'s are standard Normal distributions.
More importantly, we have cov$\{  (X_1-\mu) + (X_2-\mu),   (X_1-\mu) - (X_2-\mu)     \} = 0$, and by bivariate Normality we have $ (X_1-\mu) + (X_2-\mu) \ind  (X_1-\mu) - (X_2-\mu) $ or $Z_1\ind Z_2$ (``$\ind$'' for independence). 
The independence implies that $C = Z_1/Z_2$ follows a standard Cauchy distribution, which further leads us from
Equation (\ref{eq::fisher-z}) to
$$
\arctanh(R) \sim \log \Big|  \frac{Z_1 }{Z_2}  \Big| + \arctanh(\rho)  \sim \log |C| + \arctanh(\rho)\sim  \frac{\pi}{2} Y  + \arctanh(\rho) .
$$
Therefore, $\arctanh(R)$ is distributed around $\arctanh(\rho)$ with variability induced by a $\pi Y/2$ random variable.
Define $ V =\arctanh(R) $ and $\xi = \arctanh(\rho)$. Applying the density formula for the location and scale transformation (Casella and Berger 2001, pp 116), we obtain the density of $V$:
$$
f_V(v) =      \frac{1}{\pi}  \text{sech} (v - \xi), \quad v\in (-\infty, +\infty).
$$

Historically, Fisher (1921) used a geometrical approach to obtain the density function of the ``intraclass correlation coefficient'' for $n$ iid pairs of bivariate Normal variables. Thorndike (1905)'s measure of similarity between twins is a special case with $n=1$, and Fisher (1921) obtained the density of $\arctanh(R)$ as an application of his result.
Here, we offer a new elementary but more transparent proof of the density of $\arctanh(R)$.

\section{How Informative is the Jeffreys' Prior for Contingency Tables?}

\begin{table}[ht]
\caption{A $2\times 2$ Contingency Table}\label{tb::2x2}
\centering
\begin{tabular}{|c|cc|c|}
\hline
  &  $D=1$  & $D=0$ &  row sum \\
\hline
$E=1$ & $n_{11} $  & $n_{10}$ & $n_{1+}$\\
$E=0$ & $n_{01} $  & $n_{00} $ & $n_{0+}$\\
\hline
column sum & $n_{+1}$ & $n_{+0}$ & $n_{++}$  \\
\hline
\end{tabular}
\end{table}

Two by two contingency tables have very wide applications. For example, epidemiologists are interested in the association between a binary exposure $E$ and a binary outcome $D$.
As shown in Table \ref{tb::2x2}, the cell probabilities $p_{ij} = P(E=i,D=j)$ are of primary interest based on the observations of the cell counts $n_{ij} = \#\{k:E_k= i, D_k = j \}$. We first assume that $(n_{11}, n_{10}, n_{01}, n_{00})$ follows a Multinomial distribution with parameter $\bm{p} = (p_{11}, p_{10}, p_{01}, p_{00})$.
One commonly-used ``non-informative'' prior satisfying invariance under reparametrization is the Jeffreys' prior, which is proportional to the square root of the determinant of the Fisher information $\sqrt{\text{det}\{ I(\bm{p})\} } $.
For the Multinomial model, we can verify that
the Jeffreys' prior for $\bm{p} $ is
Dirichlet$(1/2,1/2,1/2,1/2)$.
A detailed discussion of Jeffreys' prior and the Multinomial model can be found in Section 1.3 of Box and Tiao (1973).

It is relatively direct to obtain that the prior of the marginal probability of $E$ is $p_{1+} = p_{11} + p_{10}\sim $ Beta$(1,1)\sim $ Uniform$(0, 1)$, and the prior of the marginal probability of $D$ is also $p_{+1} = p_{11} + p_{01}\sim$ Beta$(1,1)\sim $ Uniform$(0,1)$.
Since the priors of the marginal distributions are uniform,
one natural question is about the implied prior distribution of the log odds ratio defined as $W \equiv \log \{p_{11}p_{00}/ (p_{10}p_{01}) \}$.
How informative is the Jeffreys' prior for $W$?

The following representation of the
Dirichlet distribution is crucial for our proof:
\begin{eqnarray}
\label{eq::diri-gamma}
(p_{11}, p_{10}, p_{01}, p_{00})\sim \frac{(X_{11}, X_{10}, X_{01}, X_{00})}{X_{11}  + X_{10} + X_{01} + X_{00}},
\end{eqnarray}
where $X_{ij} \stackrel{iid}{\sim} \Gammad(1/2, 1) \sim \Gammad(1/2, 1/2)/2$, and the last equation is due to the scale transformation of the Gamma distribution.
Since $\chi^2_1$ distribution is a special Gamma distribution with $\chi^2_1\sim $ Gamma$(1/2,1/2)$ (Casella and Berger 2001, pp 101), we have $X_{ij}\stackrel{iid}{\sim} Z_{ij}^2/2$, where $Z_{ij}\stackrel{iid}{\sim} \mathcal{N}(0,1)$.
With these ingredients, we have the following ``one-line'' proof for the distribution of $W$:
\begin{eqnarray}
\label{eq::rep}
W \sim \log \frac{X_{11} X_{00} }{X_{10} X_{01} } 
\sim  \log \frac{Z_{11}^2 Z_{00}^2}{Z_{10}^2 Z_{01}^2} 
=2 \log\Big| \frac{Z_{11} }{ Z_{10} } \Big| + 2\log \Big|\frac{Z_{00} }{Z_{01}} \Big|
\sim  2 \log |C_1| + 2\log |C_2| 
\sim  \pi (Y_1 + Y_2),
\end{eqnarray}
where $Z_{11}/Z_{10}\sim C_1$ and $ Z_{00}/Z_{01}\sim C_2$ are iid standard Cauchy distributions, and $2 \log|C_i| / \pi  \sim Y_i $ are iid HS distributions due to the representation (\ref{eq::HS-Cauchy}).
According to (\ref{eq::rep}), we can see that 
$W$ is symmetric with mean $0$ and variance $2\pi^2 \approx 19.74$.
Baten (1934) derived the density function of the sum of $n$ iid HS random variables. Applying his result to (\ref{eq::rep}) with $n=2$, 
we can obtain the density of $W$:
\begin{eqnarray*}
\label{eq::dense}
f_W(w) = \frac{w}{2\pi^2}  \text{csch}\left( \frac{w}{2} \right) = \frac{w}{\pi^2 (e^{w/2} - e^{-w/2})}, \quad  w\in (-\infty, +\infty).
\end{eqnarray*}

The discussion above holds also for Binomial sampling with $n_{11}\sim $ Binomial$(n_{1+}, q_1)$, $n_{01}\sim $ Binomial$(n_{0+}, q_0)$, and $n_{11}\ind n_{01}$, where the marginal counts $(n_{1+}, n_{0+})$ are fixed and probabilities $(q_1, q_0)$ are unknown parameters. The Jeffreys' priors are $q_1\sim $ Beta$(1/2,1/2)$, $q_0\sim $ Beta$(1/2,1/2)$, and they are independent. Similar to (\ref{eq::diri-gamma}), Beta distributions can also be represented by Gamma distributions, e.g., $q_1\sim X_{11}/(X_{11}+X_{10})$ and $q_0\sim X_{01}/(X_{01}+X_{00})$, where $X_{ij}\stackrel{iid}{\sim} \Gammad(1/2,1/2)$ as defined before.
Consequently, the Jeffreys' prior for the log odds ratio is $ \log [ q_1 (1-q_0) / \{  q_0(1-q_1)  \} ] \sim \log \{  X_{11}X_{00} / ( X_{10}X_{01} ) \}$, which follows the same distribution as (\ref{eq::rep}).

%

\section{What If the Instrumental Variable is Invalid?}

Causal inference from observational studies often suffers from selection bias, but randomized experiment may not be feasible due to ethical or logistic problems. Encouragement experiments are attractive tools, when direct manipulation of the treatment is impossible but encouragement of treatment is feasible. 
Angrist, Imbens and Rubin (1996) provided a formal causal framework using potential outcomes, and we will briefly review their main results.
For example, in order to evaluate the causal effect of a job training program on the log of the wage, defined as $Y_i$, we randomly encourage individuals to take the job training program ($T_i=1$ if individual $i$ is encouraged and $0$ otherwise). However, the treatment received by individual $i$, $D_i$, may be different from $T_i$ ($D_i=1$ is individual $i$ receives the job training program and $0$ otherwise). 
Let $\{D_i(1), D_i(0)\}$ and $\{Y_i(1), Y_i(0)\}$
be the potential outcomes of the treatment received and log wage with and without encouragement. 
The monotonicity assumption requires $D_i(1)\geq D_i(0)$ for all $i$, which implies that encouragement will not make each individual $i$ less likely to take the treatment. And they further assumed the ``exclusion restriction'' (ER): when $D_i(1)=D_i(0)$, we have $Y_i(1) = Y_i(0)$. 
ER implies that the encouragement changes the log wage for individual $i$ only when the encouragement changes his/her potential outcomes of treatment received.
When ER holds, $T$ is called an instrumental variable (IV).
Under the assumptions of randomization of $T$, monotonicity and ER, Angrist, Imbens and Rubin (1996) showed that the traditional IV estimator
$$
\widehat{\beta}_{IV} = \frac{\overline{Y}_1 - \overline{Y}_0 }{ \overline{D}_1 - \overline{D}_0}
$$
has a valid causal interpretation,
where $\overline{Y}_1$ is the sample mean of $Y$ under treatment, and other quantities are defined analogously.
It consistently estimates the complier average causal effect:
$$
CACE \equiv E\{ Y_i(1) - Y_i(0)\mid D_i(1) = 1, D_i(0)=0\},
$$
which is the average causal effect of $T$ on $Y$ for the ``compliers'', i.e., the individuals with $D_i(1) = 1$ and $ D_i(0)=0$. Since $T_i=D_i$
for all compliers, $CACE$ can also be interpreted as the average causal effect of $D$ on $Y$.
Under regularity conditions, the IV estimator is consistent for $CACE$ and asymptotically Normal.
However, in many practical problems, the IV is very ``weak'', in the sense that $  \overline{D}_1 - \overline{D}_0 $ is very close to zero. In this case, researchers doubt the validity of the IV, since it is possible that $ T\ind D $ and the small realized value of $ \overline{D}_1 - \overline{D}_0$ is only noise.
In the following, we will obtain the asymptotic distribution of the IV estimator, if the IV is invalid, i.e., $T\ind D$.

As the sample size $N\rightarrow \infty$, we have
\begin{eqnarray*}
\sqrt{N}
\begin{pmatrix}
\overline{Y}_1 - \overline{Y}_0\\
\overline{D}_1 - \overline{D}_0
\end{pmatrix}
\stackrel{d}{\longrightarrow} \mathcal{N}_2
\left\{ 
\begin{pmatrix}
0\\
0
\end{pmatrix},
\frac{1}{\pi_1 \pi_0}
\begin{pmatrix}
\sigma_Y^2  &  \sigma_{YD} \\
\sigma_{YD}  & \sigma^2_D
\end{pmatrix}
\right\},
\end{eqnarray*}
where $\pi_1$ and $\pi_0$ are probabilities of getting treatment and control, and $(\sigma_Y^2, \sigma_D^2, \sigma_{YD})$ are variances and covariance of $Y$ and $D$.
To simplify the derivation, we define $\rho_{YD} = \sigma_{YD}/(\sigma_Y \sigma_D)$ as the correlation coefficient between $Y$ and $D$, and let $Z_1$ and $Z_2$ be two independent $\mathcal{N}(0,1)$.
We can represent the above asymptotic distributions as $ \sqrt{N\pi_1\pi_0} ( \overline{D}_1 - \overline{D}_0)  \stackrel{d}{\rightarrow} \sigma_D Z_2$ and $   \sqrt{N\pi_1\pi_0}  (\overline{Y}_1 - \overline{Y}_0)  \stackrel{d}{\rightarrow} \sigma_Y \left(  \rho_{YD} Z_2 + \sqrt{1 - \rho_{YD}^2} Z_1 \right) $, which characterize both the marginal and joint asymptotic distributions.
Using the continuous mapping theorem, we have
\begin{eqnarray}
\label{eq::iv}
\widehat{\beta}_{IV}  
= \frac{\sqrt{N\pi_1\pi_0}  (\overline{Y}_1 - \overline{Y}_0) }{\sqrt{N\pi_1\pi_0} ( \overline{D}_1 - \overline{D}_0) } 
\stackrel{d}{\longrightarrow} \frac{   \sigma_Y   \left( \rho_{YD}  Z_2 + \sqrt{1 - \rho_{YD}^2   } Z_1 \right)   }{  \sigma_D  Z_2  } 
\sim  \rho_{YD}\frac{\sigma_Y}{\sigma_D} + \eta \frac{Z_1}{Z_2} 
\sim  \rho_{YD}\frac{\sigma_Y}{\sigma_D}  + \eta  C,
\end{eqnarray} 
where $ C = Z_1/Z_2\sim $ Standard Cauchy, and $\eta = \sigma_Y\sqrt{1-\rho_{YD}^2}/\sigma_D$.
From (\ref{eq::iv}), the asymptotic mean of $\widehat{\beta}_{IV}$ is infinity, although it is centered at $\rho_{YD}\sigma_Y / \sigma_D$, i.e., the probability limit of the least square (LS) estimator of $Y$ on $D$, $\widehat{\beta}_{LS}$.
Therefore, the IV estimator can be even worse than the LS estimator in presence of an invalid IV.
Under mild regularity conditions, $\widehat{\beta}_{LS}$ is consistent for $\rho_{YD} \sigma_Y / \sigma_D$, i.e.,
$\widehat{\beta}_{LS} - \rho_{YD} \sigma_Y / \sigma_D = o_P(1)$, and therefore we have
$
\widehat{\beta}_{IV} - \widehat{\beta}_{LS} \stackrel{d}{\longrightarrow} \eta C
$
and 
\begin{eqnarray*}
\label{eq::asym}
\log |\widehat{\beta}_{IV} - \widehat{\beta}_{LS} | \stackrel{d}{\longrightarrow} \log \eta  + \frac{\pi}{2}  Y,
\end{eqnarray*}
where $Y\sim HS.$
Putting the difference of the LS estimator and IV estimator on the log scale, it is a location and scale transformation of the HS distribution, with mean $\log \eta$ and variance $ \pi^2/4.$

\section*{References}
\begin{description} \itemsep=-\parsep \itemindent=-1.2 cm

\item Angrist, J. D., Imbens, G. W., and Rubin, D. B. (1996). Identification of Causal Effects Using Instrumental Variables. {\it Journal of the American statistical Association}, {\bfseries 91}, 444-455.

\item Baten, W. D. (1934). The Probability Law for the Sum of $n$ Independent Variables, Each Subject to the Law $(2h)^{-1}\text{sech}(\pi x/2h)$. {\it Bulletin of the American Mathematical Society}, {\bfseries 40}, 284-290.

\item Box, G. E. P. and Tiao, G. C. (1973). {\it Bayesian inference in statistical analysis.} John Wiley, New York.

\item Casella, G. and Berger, R. (2001). {\it Statistical Inference}, 2nd edn. Brooks/Cole.

\item Fisher, R. A. (1921). One the Probable Error of a Coefficient of Correlation Deduced from a Small Sample. {\it Metron}, {\bfseries 1}, 1-32.

\item Harkness, W. L., and Harkness, M. L. (1968). Generalized Hyperbolic Secant Distributions. {\it Journal of the American Statistical Association}, {\bfseries 63}, 329-337.

\item Manoukian, E. B. and Nadeau, P. (1988). A Note on the Hyperbolic-Secant Distribution. {\it The American Statistician}, {\bfseries 42}, 77-79.

\item Morris, C. N. (1982). Natural Exponential Families with Quadratic Variance Functions. {\it The Annals of Statistics}, {\bfseries 10}, 65-80.

\item Morris, C. N. and Lock, K. F. (2009). Unifying the Named Natural Exponential Families and their Relatives. {\it The American Statistician}, {\bfseries 63}, 247-253.

\item Perks, W. F. (1932). On Some Experiments in the Graduation of Mortality Statistics. {\it Journal of the Institute of Actuaries}, {\bfseries 63}, 12-57.

\item Talacko, J. (1956). Perks' Distributions and Their Role in the Theory of Wiener's Stochastic Variables. {\it Trabajos de Estadistica}, {\bfseries 17}, 159-174.

\item
Thorndike, E. L. (1905). Measurement of Twins. {\it The Journal of Philosophy, Psychology and Scientific Methods}, {\bfseries 2}, 547-553.

\item Vaughan, D. C. (2002). The Generalized Secant Hyperbolic Distribution and Its 
Properties. {\it Communications in Statistics-Theory and Methods}, {\bfseries 31}, 219-238.

\end{description}

\end{document}